\documentclass[12pt]{article}

\usepackage{amsthm}
\usepackage{amsmath}
\usepackage{amssymb}

\newcommand{\R}{{\mathbb R}}

\newcommand{\Z}{{\mathbb Z}}
\newcommand{\N}{{\mathbb N}}
\newcommand{\X}[2]{\ensuremath{\hat\chi^{(#1)}(#2)}}
\newcommand{\Xone}[1]{\ensuremath{\hat\chi(#1)}}

\newcommand{\colors}{\ensuremath{C_1,\ldots,C_m}}
\newcommand{\by}{\times}

\newtheorem*{thm}{Theorem}
\newtheorem*{prop}{Proposition}
\newtheorem{prob}{Problem}

\theoremstyle{definition}
\newtheorem*{defns}{Definitions}

\begin{document}
\title{Upper Chromatic Numbers:  An Update\footnote{Final version published in \emph{Geombinatorics,} vol.~10 no.~1 (2000), pp.~4--11.}}
\author{Aaron Abrams \thanks{I thank Pete Johnson for his encouragement, assistance, and patience
throughout the writing of this paper.}
\\
Department of Mathematics\\
University of California, Berkeley \\
abrams@math.berkeley.edu
}
\date{}
\maketitle

\section{Introduction}

In an earlier volume of this journal, Greenwell and Johnson (\cite{g&j})
defined a new chromatic quantity which they called ``strange but inevitable."
Since then, a small amount of progress has been made in understanding
these ``upper" chromatic numbers.  This paper attempts to describe the
current state of knowledge in this coloring topic.

I suspect that most readers of this journal are familiar with chromatic
numbers, both for graphs and for subsets of Euclidean space.  But to
establish some context, and for those who aren't, I'll review the 
definitions.  For a graph $G$, the {\em chromatic number} $\chi(G)$
is the minimum number of colors necessary to color the vertices of $G$
such that no two adjacent vertices receive the same color.  For a set
$S\subseteq \R^n$ (or indeed any metric space), the {\em chromatic number}
$\chi(S)$ of $S$ is the minimum number of colors necessary to color
the points of $S$ such that no two points of $S$ which are a distance 
1 from each other are assigned the same color.

%
%
[Note that each of the above chromatic numbers can be viewed as 
a special case of the other; for $\chi(S)$ is the chromatic number
of the {\em unit distance graph} of $S$, and conversely $\chi(G)
=\chi(V(G))$ where $G$ is thought of as a metric space in which 
each edge has length 1.  Here $V(G)$ denotes the set of vertices of $G$.]

For future reference, we remark that $\chi(G)$
is obviously finite whenever $G$ is finite; $\chi(G)$ is also finite
if the degrees of the vertices are bounded.  Similarly, $\chi(S)$
is finite for any $S\subseteq\R^n$.  However for most choices of $S$,
the exact value of $\chi(S)$ is unknown.  In particular $\chi(\R)=2$
but the famous and elementary result of Hadwiger, Isbell, and Nelson (see 
\cite{soifer} and \cite{hadwiger} or \cite{had2})
that $4\leq\chi(\R^2)\leq 7$ provides the best known bounds in two
dimensions.

\section{Upper chromatic numbers}

From now on we consider subsets $S$ of $\R^n$.
Looking again at the definition of $\chi(S)$,
our requirement for colorings of $S$ was that no two points
a unit distance apart be colored with the same color.  Said differently,
we require that no two points of the same color be a distance 1 apart.  
This slight change of perspective invites us to generalize our acceptable
colorings
by allowing different distances to be ``forbidden" for different colors.  
For instance, what if we say that
no two red points should be distance 1 apart but no two blue points
should be distance 2 apart?  Can we color $S$ red and blue subject to
these constraints?  

This type of question was discussed by Greenwell and Johnson
in \cite{g&j}.  They note that
although $\chi(\R)=2$, $\R$ cannot be colored red and blue subject
to the constraints of the previous paragraph.  (Why not?)  This leads them to
ask the following:  How many colors are needed such that {\em whatever}
distance is forbidden for each color, $S$ can be colored subject
to the given constraints?  Thus the game is as follows:  
first you choose a number $m$ of colors you want to use.
Then I restrict a distance for each color. Finally, you try to color $S$,
obeying my restrictions.  If you have chosen $m$ in such a way that 
you can succeed no matter what I say, then
the {\em upper chromatic number} $\Xone S$ is at most $m$.  
The {\em $k^{\text th}$ upper chromatic number} of $S$, written
\X k S, is defined by the same sort of game, except that I may
restrict as many as $k$ distances for each color, rather than just
one.  

To be precise, and to establish notation, we record the

\begin{defns}

Let $S\subseteq \R^n$, and let $k,m\in\N$.  A $k\by m$ array
$D=(d_{ij})$ of positive real numbers will be called a {\em 
restriction array}.  Given $D$, a coloring of $S$ with the colors
\colors\ is a {\em $D$-coloring} if for each $i$ and $j$
($1\leq i\leq k$, $1\leq j\leq m$), no two points colored $C_j$
are a distance $d_{ij}$ apart.

The $k^{\text{th}}$ {\em upper chromatic number} of $S$, written 
$\X k S$, is the smallest integer $m$ such that for any $k\by m$
restriction array $D$, there exists a $D$-coloring of $S$.
If no such integer exists then $\X k S =\infty$.  The {\em upper
chromatic number} of $S$, written \Xone S, is defined by $\Xone S = 
\X 1 S$.
%
%
%
\end{defns}

[A remark about terminology:  Johnson suggests that the {\em lower 
chromatic number} $\check\chi(S)$ of $S\subseteq \R^n$ should be the 
smallest number
$m$ such that $S$ is $D$-colorable for {\em some} $1\by m$ array $D$.
Obviously $\check\chi (S)\leq\chi(S)\leq \Xone S$.  See Soifer
(\cite{soifer1}) or any of the papers \cite{soiferetc}, \cite{raiskii}, 
\cite{woodall}, for proofs that $\check\chi(\R^2) \leq 6$.]

\subsection{Finiteness}
A subtlety of the upper chromatic number is that it is not at all apparent
that it should be finite.  For the ordinary chromatic number $\chi(S)$,
one can establish its finiteness by constructing a single coloring; 
for \X k S one needs a $D$-coloring for each restriction array D.

In \cite{g&j} Greenwell and Johnson actually construct $D$-colorings of
$\R$ for every $1\by 3$ restriction array $D$; this shows that $\Xone \R
=3$.  (Why isn't \Xone \R\ less than 3?)  This remains the only value
of $\X k {\R^n}$ which is known exactly.  This shouldn't be too 
surprising, since again even $\chi(\R^2)$ is not known exactly.
What is perhaps surprising is that for $n\geq 2$ and any $k$, we
do not even know if $\X k {\R^n}$ is finite!

In \cite{me} I showed that \X k \R\ is finite for all $k$.  My
proof relies on probabilistic methods, and is therefore usually 
not constructive.  I obtained a bound on \X k \R\
which is asymptotically $c^k k!$; in particular the best I
could do when $k=2$ was to get $\X 2 \R < 110$.

More recently Archer has shown that \X k \R\ is $O(k)$.  In his paper
\cite{archer} he establishes the upper bound $\X k \R \leq \lceil 4ek \rceil 
(\leq 11k)$, which is currently the best known bound.  This result
also lends importance to the observation made in \cite{g&j} that
$\X k \R > k$ (see below).  Archer's methods are also probabilistic.

%

The higher dimensional cases are all wide open.  Archer (\cite{archer}) 
proved inductively that \X k {\R^n} is finite if $\R^n$ is given the
$l^{\infty}$ (``sup") norm or the $l^1$ (``taxicab") norm, but
in Euclidean space even \Xone {\R^2} is not known to be finite.

\begin{prob}\label{main}  (See \cite{g&j}.)  Is $\Xone {\R^2}$ 
finite?  What about $\X k {\R^n}, n \geq 2$?
\end{prob}

\subsection{Lower bounds}

Looking for lower bounds is entertaining.  Observe that generally
$\X k {\R^n} \geq b$ if there exists a set $S$ in $\R^n$ with $|S|=b$
and $$\bigl|\{d(x,y):x,y\in S\}\bigr|\leq k.$$  Such an $S$ is called 
a {\em $k$-distance set}.  Here are just a few observations:

\begin{itemize}
\item $\X k \R > k$ (\cite{g&j}).  Proof:  The set $\{0,1,\ldots,k\}$ 
is a $k$-distance set in $\R$.
\item $\X k {\R^2} > 2k$.  Proof:  The vertices of a regular $(2k+1)$-gon
form a $k$-distance set in $\R^2$.
\item $\X k {\R^k} \geq 2^k$.  Proof:  The vertices of the unit hypercube 
$\{0,1\}^k$ form a $k$-distance set in $\R^k$.
\item $\X k {\R^n} \geq \binom{n+1}{k}$.  Proof:  It is a (challenging)
exercise to construct a $k$-distance set in $\R^n$ of size $\binom{n+1}{k}$.
(See \cite{babai}.)
\item $\X 3 {\R^3} \geq 12$.  Proof:  The vertices of an icosahedron
form a 3-distance set in $\R^3$.
\end{itemize}
I expect that bounds obtained by examining $k$-distance sets are 
in general quite bad; these sets obstruct only a few restriction arrays.
Here are two more bounds, obtained by playing around:
\begin{itemize}
\item $\X 2 \Z \geq 4$.  Proof:  $\Z$ is not $D$-colorable for $D=
\binom{1\,1\,1}{2\,3\,4}$.
\item $\Xone {\Z^2} \geq 4$.  Proof:  $\Z^2$ is not 
$(1\,\,\sqrt{2}\,\ 2)$-colorable.
\end{itemize}

\begin{prob} Improve any of these bounds, especially $k<\X k \R \leq 
\lceil 4ek \rceil$.
In particular, improve the bounds $4\leq\X 2 \R \leq 22$.  (See Problem
1 of \cite{archer}.  Also, it is known that \X k \R = \X k \Z; see below.)
\end{prob}

\section{Some variants}

\subsection{One more chromatic number}  For $S\subseteq \R^n$, and 
$k\in\N$, let $c_k(S)$ be the minimum number of colors such that
for any $k$ forbidden distances $d_1,\ldots,d_k$, $S$ can be colored
such that no two points the same color are $d_i$ apart, for any $i$.
(This definition is due to Babai, \cite{laci}).  
This ``symmetric"
version of $\X k S$ amounts to requiring that $S$ be $D$-colorable
for $k\times m$ arrays whose columns are all identical; thus
we have $\chi(S)\leq c_k(S) \leq \X k S$.  Johnson (\cite{johnson})
argues that $c_k(S)$ is finite for all $k$ and $S$; we leave
it as an exercise to prove the

\begin{thm}  If $S\subseteq \R^n$ then $c_k(S)\leq 
\left(\chi(\R^n)\right)^k$.
\end{thm}

Observe also that any lower bound on $\X k S$ established using 
$k$-distance sets will apply equally well to $c_k(S)$; this is
why I think such bounds on $\X k S$ are poor.

Of course, we don't know very much about $c_k$.  A future paper will 
discuss some results in this direction, but as usual the questions far
outnumber the answers.

\subsection{Final ruminations}
In $\R$, there is a universal upper bound (namely, 2) on the number of
points which can be a given distance from a given point.  One of the
main tools of \cite{me} and \cite{archer} is the Lov\'asz
Local Lemma (see \cite{lll1} or \cite{lll2}), which is designed for just 
such ``locally bounded" situations.
The failure of our methods to generalize
to higher dimensions can be traced to this basic principle, one form
of which is
made explicit as Corollary 5 of \cite{archer}.

Observe, for instance, that if we choose to restrict {\em vectors} 
instead of distances,
then the corresponding chromatic number (which you might write as 
\X {\vec{k}} {\R^n}) is finite.  In fact a slick argument based on
an idea of Eric Wepsic \cite{wepsic} shows that
$\X {\vec{k}} {\R^n}=\X k \R$.
(Thus \X {\vec{k}} {\R^n} doesn't depend on $n$!)  

Another application of this principle will follow shortly.  But first, 
we need one more bit of notation from \cite{g&j}.  

Let $R\subseteq(0,\infty)$.
Let $\hat\chi_R^{\phantom{()}}$ and $\hat\chi_R^{(k)}$ be defined 
as before, except
the forbidden distances should come from $R$; $R$ is the set of
``possible" restrictions.  Greenwell and Johnson proved in \cite{g&j}
that $\hat\chi_R^{\phantom{()}}(\R^n)$ is finite if $R=[\epsilon,M]$ 
for some $0<\epsilon<M<\infty$.

We now consider $\Z^n$.  Tantalizingly, Archer has proved 
that $\X k {\Z^n} 
= \X k {\R^n}$ for the $l^1$ and $l^{\infty}$ metrics on
$\R^n$; this suggests that results about $\Z^n$ may well bear on
the general case.  However, we do not know if this theorem holds
in Euclidean space as well.  

\begin{prob}  Is $\X k {\R^n}=\X k {\Z^n}$ for all $k,n$?
\end{prob}

Anyhow, the only distances
occurring in $\Z^n$ are square roots of integers; the number of
points in $\Z^n$ which are a distance $\sqrt{d}$ from a given point
is the same as the number of ways of writing $d$ as a sum of
(at most) $n$ squares.  Thus if we let
\begin{eqnarray*}
R_N=\{\sqrt{d}:&d\in\N \ \ \text{ can be written as a sum of } n 
\text{ squares}\\
&\phantom{blah } \text{ in at most } N \text{ ways}\}, 
\end{eqnarray*}
then again by our finiteness principle,
$\hat\chi^{(k)}_{R_N}(S)$ is finite for any $S\subseteq\Z^n$.
This may sound like a bizarre condition, but consider what it means when
$n=2$ and $N=1000$, say.  It seems like ``most"
distances are allowed in $D$ (though infinitely many aren't); yet
we can't prove that \Xone {\Z^2} is finite.  [Note that the set
$R_N$ depends on $n$, though it is left out of the notation.]

When one plays around a little with $\X 2 \Z$, one can easily discover
many congruence conditions on the entries of $D$ which imply that
$\Z$ is $D$-colorable.  
%
%
Here is a corresponding observation regarding \Xone {\Z^2}.  
As before, we may assume that 
$D$ is of the form $(\sqrt{d_1}, \ldots, \sqrt{d_m})$, where
the $d_i$ are natural numbers.
Consider the following property that $D$ might have: 
\begin{description}
\item[$(*)$]
if $i\not=j$ then the number of times that 2 divides $d_i$
is different from the number of times that 2 divides $d_j$.
\end{description}

After some thought one can see:
\begin{prop}
\Xone {\Z^2} is finite if and only if $\Z^2$ is $D$-colorable
for every $D$ with property $(*)$.
\end{prop}

Thorough treatments of this and related results are postponed to
a future paper.

\bibliography{upper}
\bibliographystyle{amsplain}

\end{document}